\documentclass[12pt]{article}
\usepackage{amssymb,amsmath}
\usepackage{cases}
\usepackage{amsfonts}
\usepackage{color,xcolor}
\usepackage[left=2.0cm,right=2.0cm,top=2.0cm,bottom=2.0cm]{geometry}
\usepackage[colorlinks,citecolor=blue,urlcolor=blue]{hyperref}
\usepackage[pagewise]{lineno}\modulolinenumbers

\newtheorem{theorem}{Theorem}[section]

\newtheorem{lemma}{Lemma}[section]
\newtheorem{proposition}{Proposition}[section]

\newtheorem{definition}{Definition}[section]
\newtheorem{remark}{Remark}[section]

\newcommand{\bal}{\begin{align}}
\newcommand{\bbal}{\begin{align*}}
\newcommand{\beq}{\begin{equation}}
\newcommand{\eeq}{\end{equation}}
\newcommand{\bca}{\begin{cases}}
\newcommand{\eca}{\end{cases}}

\newcommand{\R}{\mathbb{R}}

\newcommand{\bi}{\Big}

\begin{document}

\title{Non-uniform continuity on initial data for the two-component b-family system in Besov space}

\author{Xing Wu$^{1, }$\thanks{Corresponding author. ny2008wx@163.com (Xing Wu)}\;, Cui Li$^2$,  Jie Cao$^1$\\
\small \it$^1$College of Information and Management Science,
Henan Agricultural University,\\
\small  Zhengzhou, Henan, 450002, China\\
\small \it $^2$School of Mathematics and Information Science, Henan University of Economics and Law,\\
\small Zhengzhou, Henan, 450000, China}

\date{}

\maketitle\noindent{\hrulefill}

{\bf Abstract:} In this paper, we consider the Cauchy problem of a two-component b-family system, which includes the two-component Camassa-Holm system and the two-component Degasperis-Procesi system. It is shown that the  solution map of the two-component b-family system is not uniformly continuous on the initial data in Besov spaces $B_{p, r}^{s-1}(\mathbb{R})\times B_{p, r}^s(\mathbb{R})$ with $s>\max\{1+\frac{1}{p}, \frac{3}{2}\}$, $1\leq p, r< \infty$. Our result covers and extends the previous non-uniform continuity in Sobolev spaces $H^{s-1}(\mathbb{R})\times H^s(\mathbb{R})$ for $s>\frac{5}{2}$ to Besov spaces.


\vskip0mm\noindent{\hrulefill}

\section{Introduction}\label{sec1}
In this paper, we are concerned with the Cauchy problem for the following  two-component b-family system on $\mathbb{R}$:
\begin{eqnarray}\label{eq1}
        \left\{\begin{array}{ll}
          m_t=um_x+ k_1u_xm+k_2\rho\rho_x,\\
         \rho_t=k_3(u\rho)_x,\\
          m=u-u_{xx},\\
         u(0, x)=u_0, \rho(0, x)=\rho_0. \end{array}\right.
        \end{eqnarray}
which was introduced by Guha in \cite{Guha 2007}. As shown in \cite{Guha 2007},  there are two cases about this system: (i) $k_1=b$, $k_2=2b$ and $k_3=1$; (ii) $k_1=b+1,$ $k_2=2$ and $k_3=b$ with $b\in \mathbb{R}.$

If $k_1=2$ and $k_3=1$, then system (\ref{eq1}) becomes the following two-component Camassa-Holm system
\begin{eqnarray}\label{eq2}
        \left\{\begin{array}{ll}
          m_t=um_x+ 2u_xm+\sigma\rho\rho_x,\\
         \rho_t=(u\rho)_x,\\
          m=u-u_{xx}, \end{array}\right.
        \end{eqnarray}
here $\sigma=\pm 1.$ System (\ref{eq2}) was  derived by Constantin and Ivanov \cite{Constantin 2008} in the context of shallow water theory, and then has attracted much more attention. The local well-posedness for system (\ref{eq2}) in Sobolev and Besov space were established in \cite{Escher 2007, Constantin 2008, Guan 2010, Gui 2011}. The global existence of strong solutions and wave-breaking mechanism were investigated in
\cite{Escher 2007, Guan 2010, Gui 2011, Gui 2010}, and the global weak solution
has been obtained in \cite{Guan 2011}. Moreover, when $\rho=0$, (\ref{eq2}) reduces to the classical Camassa-Holm equation modeling the unidirectional propagation of shallow water waves over a flat bottom.
The Cauchy problem of the Camassa-Holm equation was studied in the
series of papers  \cite{Constantin 1998, 1Constantin 1998, Constantin 2000, 1Constantin 2000, Danchin 2001, Li 2016, Guo 2019}. Danchin\cite{Danchin 2001, Danchin 2003} showed the local existence and uniqueness of strong solutions to Camassa-Holm equation with  initial data in $B_{p, r}^s$ for  $s>\max\{1+\frac{1}{p}, \frac{3}{2}\}$, $1\leq p\leq \infty$, $1\leq r<\infty$. The continuous dependence of the solution on the initial data has been supplemented by Li and Yin in \cite{Li 2016}. Very recently, Li, Yu and Zhu \cite{Li 2020} have sharpened the results in \cite{Li 2016}   by showing that the solution map is not uniformly continuous (see \cite{Himonas 2009} for  earlier results in Sobolev space $H^s(\mathbb{R})$ with $s>\frac{3}{2}$).

If $k_1=3$ and $k_3=2$, then system (\ref{eq1}) becomes the following two-component Degasperis-Procesi system
\begin{eqnarray}\label{eq3}
        \left\{\begin{array}{ll}
          m_t=um_x+ 3u_xm+\sigma\rho\rho_x,\\
         \rho_t=2(u\rho)_x,\\
          m=u-u_{xx}, \end{array}\right.
        \end{eqnarray}
which was first proposed in \cite{Popowicz 2006} as a natural generalization of the Degasperis-Procesi equation in the context of supersymmetry. The local well-posedness of system (\ref{eq3}) was established in \cite{Jin 2010, Escher 2011, Yan 2012}, the precise blow-up scenario and some blow-up rate of strong solutions were also presented in \cite{Jin 2010, Yan 2012}. Moreover, when $\rho=0$, (\ref{eq3}) reduces to the famous Degasperis-Procesi equation. The first result of nonuniform dependence for Degasperis-Procesi equation in the periodic case with $s\geq 2$ can be found by Christov and Hakkaev
\cite{Christov 2009}. This result was lately improved to the best by Himonas and Holliman in \cite{Himonas 2011}, where non-uniform dependence on the initial data for Degasperis-Procesi equation was proven on both the circle and line for $s>\frac{3}{2}$ using the method of approximate solutions in tandem with a twisted $L^2$-norm that is conserved by the DP equation. Recently, the non-uniform continuity has been once extended to the Besov space \cite{Wu 2021}.

For $\rho\neq 0$, and $b\in \mathbb{R}$, the Cauchy problem of system (\ref{eq1}) in Sobolev space was first established by Liu and Yin \cite{Liu 2011} for $(u_0, \rho_0)\in H^s \times H^{s-1}$ with $s\geq 2$. Later,  Lv and Wang \cite{Lv 2014} proved that the solution map is not uniformly continuous for $s>\frac 52$ by using the approximate solutions.
The local well-posedness space was further enlarged,  and established in Besov space $B_{p, r}^s\times B_{p, r}^{s-1}$ with  $s>\max\{1+\frac{1}{p}, \frac{3}{2}\}$, $1\leq p, r\leq \infty$ (however, for $r=\infty$, the continuity of the data-to-solution map is established in a
weaker topology)\cite{Zhu 2013, Yan 2020}. Some aspects concerning blow-up scenario, global solutions, persistence properties and propagation speed, see
the discussions in \cite{Liu 2011, Zong 2012}.

In the present paper, motivated by \cite{Li 2020, Wu 2021},  we aim at showing that  the solution map of (\ref{eq1}) is not uniformly continuous depending on the initial data in Besov spaces $B_{p, r}^s\times B_{p, r}^{s-1}$ with $s>\max\{1+\frac{1}{p}, \frac{3}{2}\}$, $1\leq p, r< \infty$. However, one problematic issue is that we here deal with a  coupled system with these two components of the solution in different Besov spaces. On the other hand, compared with Novikov equation with cubic nonlinearity \cite{1Li 2020}, quadratic nonlinearity for system (\ref{eq1}) weakens the attenuation we need. These two aspects make the proof of several required nonlinear estimates somewhat delicate.

For studying the non-uniform continuity of a two-component b-family system, it is more convenient to express (\ref{eq1}) in the following equivalent nonlocal form
\begin{eqnarray}\label{eq4}
        \left\{\begin{array}{ll}
         u_t-uu_x=f(u)+g(\rho),\\
        \rho_t-k_3u\rho_x=k_3\rho u_x,\\
         u(0, x)=u_0, \rho(0, x)=\rho_0, \end{array}\right.
        \end{eqnarray}
where $f(u)=f_1(u)+f_2(u)$ and
\begin{eqnarray*}
f_1(u)=\partial_x(1-\partial_x^2)^{-1}(\frac{k_1}{2}u^2), \quad f_2(u)=\partial_x(1-\partial_x^2)^{-1}(\frac{3-k_1}{2}u_x^2), \quad g(\rho)=\partial_x(1-\partial_x^2)^{-1}(\frac{k_2}{2}\rho^2).
\end{eqnarray*}

 Our main result is stated as follows.
\begin{theorem}\label{the1.1} Let $s>\max\{1+\frac{1}{p}, \frac{3}{2}\}$, $1\leq p, r < \infty$. The solution map $(u_0, \rho_0)\rightarrow (u(t), \rho(t))$ of the initial value problem (\ref{eq4}) is not
uniformly continuous from any bounded subset of  $B_{p, r}^s(\mathbb{R})\times B_{p, r}^{s-1}(\mathbb{R})$ into $\mathcal{C}([0, T];  B_{p, r}^s(\mathbb{R})\times B_{p, r}^{s-1}(\mathbb{R}))$. More precisely,
there exist two sequences of solutions $(u_n^1(t), \rho_n^1(t))$ and $(u_n^2(t), \rho_n^2(t))$ such that the corresponding initial data satisfy
\begin{eqnarray*}
        \|u_n^1(t), u_n^2(t)\|_{B_{p, r}^s}+\|\rho_n^1(t), \rho_n^2(t)\|_{B_{p, r}^{s-1}}\lesssim 1,\\
\end{eqnarray*}
       \mbox{and}
   \begin{eqnarray*}
        \lim_{n\rightarrow \infty}(\|u_n^1(0)-u_n^2(0)\|_{B_{p, r}^s}+\|\rho_n^1(0)-\rho_n^2(0)\|_{B_{p, r}^{s-1}})=0,
        \end{eqnarray*}
but
\begin{eqnarray*}
  \liminf_{n\rightarrow \infty} \|u_n^1(t)-u_n^2(t)\|_{B_{p, r}^s}\gtrsim t, \quad    \liminf_{n\rightarrow \infty} \|\rho_n^1(t)-\rho_n^2(t)\|_{B_{p, r}^{s-1}}\gtrsim t, \quad t\in [0, T_0],
        \end{eqnarray*}
with small positive time $T_0$ for $T_0\leq T$.
\end{theorem}
\begin{remark}\label{rem1}
Since $B_{2, 2}^s=H^s$, our result covers and extends the previous non-uniform continuity of solutions on initial data in Sobolev spaces $H^s(\mathbb{R})\times H^{s-1}(\mathbb{R})$ for $s>\frac{5}{2}$ \cite{Lv 2014} to Besov spaces.
\end{remark}

{\bf Notations}:  Given a Banach space $X$, we denote the norm of a function on $X$ by $\|\|_{X}$, and \begin{eqnarray*}
\|\cdot\|_{L_T^\infty(X)}=\sup_{0\leq t\leq T}\|\cdot\|_{X}.
\end{eqnarray*}
For $\mathbf{f}=(f_1, f_2,...,f_n)\in X$,
\begin{eqnarray*}
\|\mathbf{f}\|_{X}^2=\|f_1\|_{X}^2+\|f_2\|_{X}^2+...+\|f_n\|_{X}^2.
\end{eqnarray*}
The symbol
$A\lesssim B$ means that there is a uniform positive constant $C$ independent of $A$ and $B$ such that $A\leq CB$.

\section{Littlewood-Paley analysis}\label{sec2}
\setcounter{equation}{0}
In this section, we will review the definition of Littlewood-Paley decomposition and nonhomogeneous Besov space, and then list some useful properties. For more details, the readers can refer to \cite{Bahouri 2011}.

There exists a couple of smooth functions $(\chi,\varphi)$ valued in $[0,1]$, such that $\chi$ is supported in the ball $\mathcal{B}\triangleq \{\xi\in\mathbb{R}:|\xi|\leq \frac 4 3\}$, $\varphi$ is supported in the ring $\mathcal{C}\triangleq \{\xi\in\mathbb{R}:\frac 3 4\leq|\xi|\leq \frac 8 3\}$. Moreover,
$$\forall\,\, \xi\in\mathbb{R},\,\, \chi(\xi)+{\sum\limits_{j\geq0}\varphi(2^{-j}\xi)}=1,$$
$$\forall\,\, \xi\in\mathbb{R}\setminus\{0\},\,\, {\sum\limits_{j\in \mathbb{Z}}\varphi(2^{-j}\xi)}=1,$$
$$|j-j'|\geq 2\Rightarrow\textrm{Supp}\,\ \varphi(2^{-j}\cdot)\cap \textrm{Supp}\,\, \varphi(2^{-j'}\cdot)=\emptyset,$$
$$j\geq 1\Rightarrow\textrm{Supp}\,\, \chi(\cdot)\cap \textrm{Supp}\,\, \varphi(2^{-j}\cdot)=\emptyset.$$
Then, we can define the nonhomogeneous dyadic blocks $\Delta_j$ as follows:
$$\Delta_j{u}= 0,\,\, \text{if}\,\, j\leq -2,\quad
\Delta_{-1}{u}= \chi(D)u=\mathcal{F}^{-1}(\chi \mathcal{F}u),$$
$$\Delta_j{u}= \varphi(2^{-j}D)u=\mathcal{F}^{-1}(\varphi(2^{-j}\cdot)\mathcal{F}u),\,\, \text{if} \,\, j\geq 0.$$

\begin{definition}[\cite{Bahouri 2011}]\label{de2.1}
Let $s\in\mathbb{R}$ and $1\leq p,r\leq\infty$. The nonhomogeneous Besov space $B^s_{p,r}(\mathbb{R})$ consists of all tempered distribution $u$ such that
\begin{align*}
\|u\|_{B^s_{p,r}(\mathbb{R})}\triangleq \Big|\Big|(2^{js}\|\Delta_j{u}\|_{L^p(\mathbb{R})})_{j\in \mathbb{Z}}\Big|\Big|_{\ell^r(\mathbb{Z})}<\infty.
\end{align*}
\end{definition}

In the following, we list some basic lemmas and properties about Besov space which will be frequently used in proving our main result.

\begin{lemma}(\cite{Bahouri 2011})\label{lem2.1}
 (1) Algebraic properties: $\forall s>0,$ $B_{p, r}^s(\mathbb{R})$ $\cap$ $L^\infty(\mathbb{R})$ is a Banach algebra. $B_{p, r}^s(\mathbb{R})$ is a Banach algebra $\Leftrightarrow B_{p, r}^s(\mathbb{R})\hookrightarrow L^\infty(\mathbb{R})\Leftrightarrow s>\frac{1}{p}$ or $s=\frac{1}{p},$ $r=1$.\\
 (2) For any $s>0$ and $1\leq p,r\leq\infty$, there exists a positive constant $C=C(s,p,r)$ such that
$$\|uv\|_{B^s_{p,r}(\mathbb{R})}\leq C\Big(\|u\|_{L^{\infty}(\mathbb{R})}\|v\|_{B^s_{p,r}(\mathbb{R})}+\|v\|_{L^{\infty}(\mathbb{R})}\|u\|_{B^s_{p,r}(\mathbb{R})}\Big).$$
(3) Let $m\in \mathbb{R}$ and $f$ be an $S^m-$ multiplier (i.e., $f: \mathbb{R}^d\rightarrow \mathbb{R}$ is smooth and satisfies that $\forall \alpha\in \mathbb{N}^d$, there exists a constant $\mathcal{C}_\alpha$ such that $|\partial^\alpha f(\xi)|\leq \mathcal{C}_\alpha(1+|\xi|)^{m-|\alpha|}$ for all $\xi \in \mathbb{R}^d$). Then the operator $f(D)$ is continuous from $B_{p, r}^s(\mathbb{R}^d)$ to $B_{p, r}^{s-m}(\mathbb{R}^d)$.\\
(4) Let  $1\leq p, r\leq \infty$ and $s>\max\{1+\frac{1}{p}, \frac{3}{2}\}$. Then  we have
$$\|uv\|_{B_{p, r}^{s-2}(\mathbb{R})}\leq C\|u\|_{B_{p, r}^{s-2}(\mathbb{R})}\|v\|_{B_{p, r}^{s-1}(\mathbb{R})}.$$
\end{lemma}

\begin{lemma}\label{lem2.2}(\cite{Bahouri 2011, Li 2017})
Let $1\leq p,r\leq \infty$. Assume that
\begin{eqnarray*}
\sigma> -\min\{\frac{1}{p}, 1-\frac{1}{p}\} \quad \mathrm{or}\quad \sigma> -1-\min\{\frac{1}{p}, 1-\frac{1}{p}\}\quad \mathrm{if} \quad \mathrm{div\,} v=0.
\end{eqnarray*}
There exists a constant $C=C(p,r,\sigma)$ such that for any solution to the
following linear transport equation:
\begin{equation*}
\partial_t f+v\cdot\nabla f=g,\qquad
f|_{t=0} =f_0,
\end{equation*}
the following statements hold:
\begin{align*}
\sup_{s\in [0,t]}\|f(s)\|_{B^{\sigma}_{p,r}}\leq Ce^{CV_{p}(v,t)}\Big(\|f_0\|_{B^\sigma_{p,r}}
+\int^t_0\|g(\tau)\|_{B^{\sigma}_{p,r}}d \tau\Big),
\end{align*}
with
\begin{align*}
V_{p}(v,t)=
\begin{cases}
\int_0^t \|\nabla v(s)\|_{B^{\sigma-1}_{p,r}}ds, &\quad \mathrm{if} \;\sigma>1+\frac{1}{p}\ \mathrm{or}\ \{\sigma=1+\frac{1}{p} \mbox{ and } r=1\},\\
\int_0^t \|\nabla v(s)\|_{B^{\sigma}_{p,r}}ds,&\quad\mathrm{if} \; \sigma=1+\frac{1}{p} \quad \mathrm{and} \quad r>1,\\
\int_0^t \|\nabla v(s)\|_{B^{\frac{1}{p}}_{p,\infty}\cap L^\infty}ds,&\quad\mathrm{if} \; \sigma<1+\frac{1}{p}.
\end{cases}
\end{align*}
\end{lemma}

\section{Non-uniform continuous dependence}\label{sec3}
\setcounter{equation}{0}
In this section, we will give the proof of our main theorem.  For brevity, we sometimes use $u_{0, n}^i$, $\rho_{0, n}^i$ to denote $u^i_n(0)$ and $\rho^i_n(0)$ respectively, $i=1, 2.$

Let $\hat{\phi}\in \mathcal{C}^\infty_0(\mathbb{R})$ be an even, real-valued and non-negative funtion on $\R$ and satisfy
\begin{numcases}{\hat{\phi}(x)=}
1, &if $|x|\leq \frac{1}{4}$,\nonumber\\
0, &if $|x|\geq \frac{1}{2}$.\nonumber
\end{numcases}
Define the high frequency function $f_n$ and the low frequency function $g_n$ by
$$f_n=2^{-ns}\phi(x)\sin \bi(\frac{17}{12}2^nx\bi), \qquad g_n=2^{-n}\phi(x), \quad n\gg1.$$
It has been showed in \cite{Li 2020} that $\|f_n\|_{B_{p, r}^\sigma}\lesssim 2^{n(\sigma-s)}$.

 Let
$$(u_n^1(0),\; \rho_n^1(0))=(f_n, 2^nf_n), \quad (u_n^2(0),\; \rho_n^2(0))=(f_n+g_n, 2^nf_n+g_n),$$
then it is easy to verify that
\begin{align}
\|u_n^1(0), \; u_n^2(0)\|_{B_{p, r}^{s+\sigma}}\lesssim 2^{n\sigma}\; \mathrm{for}\; \;\sigma\geq -1 \;\;\;\mathrm{and} \;\;\;\|\rho_n^1(0),\; \rho_n^2(0)\|_{B_{p, r}^{s+l}}\lesssim 2^{n(l+1)}\;\mathrm{for}\; \; l\geq -2,\label{eq3.1}
\end{align}

Consider the system (\ref{eq4}) with initial data $(u_n^1(0), \rho_n^1(0))$ and $(u_n^2(0), \rho_n^2(0))$, respectively.  According to the local well-posedness result in \cite{Zhu 2013, Yan 2020}, there exists corresponding solution $(u_n^1, \rho_n^1)$, $(u_n^2, \rho_n^2)$ belonging to
$\mathcal{C}([0, T];  B_{p, r}^s\times B_{p, r}^{s-1})$ and has common lifespan $T\thickapprox 1$. Moreover, by Lemma \ref{lem2.1}-\ref{lem2.2}, there holds
 \begin{eqnarray}
     &\|u_n^1\|_{L_T^\infty(B_{p, r}^{s+k})}+\|\rho_n^1\|_{L_T^\infty(B_{p, r}^{s+k-1})}\lesssim  \|u_n^1(0)\|_{B_{p, r}^{s+k}}+\|\rho_n^1(0)\|_{B_{p, r}^{s+k-1}}\lesssim2^{nk},\quad k\geq -1,\label{eq3.2}\\
      & \|u_n^2\|_{L_T^\infty(B_{p, r}^{s+l})}+\|\rho_n^2\|_{L_T^\infty(B_{p, r}^{s+l-1})}\lesssim  \|u_n^2(0)\|_{B_{p, r}^{s+l}}+\|\rho_n^2(0)\|_{B_{p, r}^{s+l-1}}\lesssim2^{nl},\quad l\geq -1.\label{eq3.3}
        \end{eqnarray}

In the following, we shall firstly show that for the selected high  frequency initial data $(u_n^1(0), \rho_n^1(0))$, the corresponding solution $(u_n^1, \rho_n^1)$ can be approximated by the initial data. More precisely, that is
\begin{proposition}\label{pro1}
 Under the assumptions of Theorem \ref{the1.1}, we have
      \begin{eqnarray}
     \|u_n^1-u_n^1(0)\|_{L_T^\infty(B_{p, r}^s)}+ \|\rho_n^1-\rho_n^1(0)\|_{L_T^\infty(B_{p, r}^{s-1})}\lesssim 2^{-\frac{n}{2}(s-\frac 32)}.\label{eq3.4}
        \end{eqnarray}
\end{proposition}
\noindent{\bf Proof} \;
Denote $\epsilon=u_n^1-u_n^1(0)$, $\delta=\rho_n^1-\rho_n^1(0)$, then we can derive from (\ref{eq4}) that $(\epsilon, \varrho)$ satisfies
\begin{eqnarray}\label{eq3.5}
        \left\{\begin{array}{ll}
        \epsilon_t-u_n^1\partial_x\epsilon=(u_n^1-u^1_{0, n})\partial_xu^1_{0, n}
        +[f(u_n^1)-f(u^1_{0, n})]+[g(\rho_n^1)-g(\rho^1_{0, n})]\\
        \qquad\qquad\qquad+f(u^1_{0, n})+g(\rho^1_{0, n})+u^1_{0, n}\partial_xu^1_{0, n},\\
        \delta_t-k_3u_n^1\partial_x\delta=k_3(u_n^1-u^1_{0, n})\partial_x\rho^1_{0, n}
        +k_3u^1_{0, n}\partial_x\rho^1_{0, n}+k_3\rho_n^1\partial_xu_n^1,\\
          \epsilon(0, x)=0, \delta(0, x)=0, \end{array}\right.
        \end{eqnarray}

Applying Lemma \ref{lem2.2} yields
\begin{align}
     \|\epsilon\|_{B_{p, r}^{s-1}}&\lesssim \int_0^t\|\partial_xu_n^1\|_{B_{p, r}^{s-1}}\|\epsilon\|_{B_{p, r}^{s-1}}d\tau+\int_0^t \|(u_n^1-u^1_{0, n})\partial_xu^1_{0, n}\|_{B_{p, r}^{s-1}}d\tau \nonumber\\
    &\quad +\int_0^t \|f(u_n^1)-f(u^1_{0, n}), \;g(\rho_n^1)-g(\rho^1_{0, n})\|_{B_{p, r}^{s-1}}d\tau \nonumber\\
     &\quad+t\|f(u^1_{0, n}),\;g(\rho^1_{0, n}),\;u^1_{0, n}\partial_xu^1_{0, n}\|_{B_{p, r}^{s-1}},\label{eq3.6}\\
     \|\delta\|_{B_{p, r}^{s-2}}&\lesssim \int_0^t\|\partial_xu_n^1\|_{B_{p, r}^{s-1}}\|\delta\|_{B_{p, r}^{s-2}}d\tau+\int_0^t \|(u_n^1-u^1_{0, n})\partial_x\rho^1_{0, n}\|_{B_{p, r}^{s-2}}d\tau \nonumber\\
    &\quad +\int_0^t \|\rho_n^1\partial_xu_n^1-\rho^1_{0, n}\partial_xu_{0, n}^1\|_{B_{p, r}^{s-2}}d\tau +t\|u^1_{0, n}\partial_x\rho^1_{0, n}, \;\rho^1_{0, n}\partial_xu_{0, n}^1\|_{B_{p, r}^{s-2}}.\label{eq3.7}
\end{align}

Using Lemma \ref{lem2.1} and  the fact that $B_{p, r}^{s-1}(\mathbb{R})$ is a  Banach algebra when $s>\max\{1+\frac{1}{p}, \frac{3}{2}\}$, we have
\begin{align*}
       \|(u_n^1-u^1_{0, n})\partial_xu^1_{0, n}\|_{B_{p, r}^{s-1}}
      &\lesssim \|u_n^1-u^1_{0, n}\|_{B_{p, r}^{s-1}}
        \|\partial_xu^1_{0, n}\|_{B_{p, r}^{s-1}}\lesssim \|u_n^1-u^1_{0, n}\|_{B_{p, r}^{s-1}}\|u^1_{0, n}\|_{B_{p, r}^s},\\
      \|f(u_n^1)-f(u^1_{0, n})\|_{B_{p, r}^{s-1}}&\lesssim\|u_n^1-u^1_{0, n}\|_{B_{p, r}^{s-1}}\|u_n^1, \; u^1_{0, n}\|_{B_{p, r}^s},\\
       \|g(\rho_n^1)-g(\rho^1_{0, n})\|_{B_{p, r}^{s-1}}&\lesssim\|\rho_n^1-\rho^1_{0, n}\|_{B_{p, r}^{s-2}}\|\rho_n^1, \;\rho^1_{0, n}\|_{B_{p, r}^{s-1}},\\
     \|(u_n^1-u^1_{0, n})\partial_x\rho^1_{0, n}\|_{B_{p, r}^{s-2}} &\lesssim\|\partial_x\rho^1_{0, n}\|_{B_{p, r}^{s-2}}\|u_n-u_{0, n}\|_{B_{p, r}^{s-1}}\lesssim \|u_n-u_{0, n}\|_{B_{p, r}^{s-1}}\|\rho^1_{0, n}\|_{B_{p, r}^{s-1}},\\
        \|u^1_{0, n}\partial_x\rho^1_{0, n}\|_{B_{p, r}^{s-2}} &\lesssim\|u^1_{0, n}\partial_x\rho^1_{0, n}\|_{B_{p, r}^{s-\frac 32}}\lesssim\|u^1_{0, n}\|_{L^\infty}\|\rho^1_{0, n}\|_{B_{p, r}^{s-\frac 12}}+\|u^1_{0, n}\|_{B_{p, r}^{s-\frac 32}}\|\partial_x\rho^1_{0, n}\|_{L^\infty},\\
        &\lesssim 2^{-ns}2^{n(-\frac 12 +1)}+2^{-\frac 32n}2^n2^{-ns}2^n\lesssim 2^{-n(s-\frac 12)},\\
      \|\rho^1_{0, n}\partial_xu_{0, n}^1\|_{B_{p, r}^{s-2}}&\lesssim\|\rho^1_{0, n}\partial_xu_{0, n}^1\|_{B_{p, r}^{s-\frac 32}}\lesssim2^{n(s-\frac 32)}\|\rho^1_{0, n}, \partial_xu_{0, n}^1\|_{L^\infty}\|\rho^1_{0, n}, \partial_xu_{0, n}^1\|_{L^p}\lesssim 2^{-n(s-\frac 12)}.
       \end{align*}
       Again using Lemma \ref{lem2.1} and the Banach algebra property of $B_{p, r}^{s-1}$, one has
       \begin{align*}
     \|f_1(u_{0, n}^1)\|_{B_{p, r}^{s-1}}&\lesssim\|(u_{0, n}^1)^2\|_{B_{p, r}^{s-2}}\lesssim\|(u_{0, n}^1)^2\|_{B_{p, r}^{s-\frac 32}}\lesssim \|u_{0, n}^1\|_{L^\infty}\|u_{0, n}^1\|_{B_{p, r}^{s-\frac 32}}\lesssim 2^{-n(s+\frac 32)},\\
     \|f_2(u_{0, n}^1)\|_{B_{p, r}^{s-1}}&\lesssim\|(\partial_xu_{0, n}^1)^2\|_{B_{p, r}^{s-2}}\lesssim \|\partial_xu_{0, n}^1\|_{L^\infty}\|\partial_xu_{0, n}^1\|_{B_{p, r}^{s-\frac 32}}\lesssim 2^{-n(s-\frac 12)},\\
     \|g(\rho_{0, n}^1)\|_{B_{p, r}^{s-1}}&\lesssim\|(\rho_{0, n}^1)^2\|_{B_{p, r}^{s-2}}\lesssim\|(\rho_{0, n}^1)^2\|_{B_{p, r}^{s-\frac 32}}\lesssim \|\rho_{0, n}^1\|_{L^\infty}\|\rho_{0, n}^1\|_{B_{p, r}^{s-\frac 32}}\lesssim 2^{-n(s-\frac 12)},\\
     \|u_{0, n}^1\partial_xu_{0, n}^1\|_{B_{p, r}^{s-1}}
      &\lesssim  \|u_{0, n}^1\|_{L^\infty} \|\partial_xu_{0, n}^1\|_{B_{p, r}^{s-1}}+\|u_{0, n}^1\|_{B_{p, r}^{s-1}}\|\partial_xu_{0, n}^1\|_{L^\infty}\lesssim 2^{-ns}.
        \end{align*}
For the term
\begin{eqnarray*}
\rho_n^1\partial_xu_n^1-\rho^1_{0, n}\partial_xu_{0, n}^1=(\rho_n^1-\rho^1_{0, n})\partial_xu_n^1
+ \rho^1_{0, n}\partial_x(u_n^1-u^1_{0, n}).
\end{eqnarray*}
Following the same procedure of estimates as above, we find that
\begin{align*}
\|(\rho_n^1-\rho^1_{0, n})\partial_xu_n^1\|_{B_{p, r}^{s-2}}&\lesssim\|\rho_n^1-\rho^1_{0, n}\|_{B_{p, r}^{s-2}}\|\partial_xu_n^1\|_{B_{p, r}^{s-1}}\lesssim\|\rho_n^1-\rho^1_{0, n}\|_{B_{p, r}^{s-2}}\|u_n^1\|_{B_{p, r}^s},\\
\|\rho^1_{0, n}\partial_x(u_n^1-u^1_{0, n})\|_{B_{p, r}^{s-2}}&\lesssim\|\partial_x(u_n^1-u^1_{0, n})\|_{B_{p, r}^{s-2}}\|\rho^1_{0, n}\|_{B_{p, r}^{s-1}}\lesssim\|u_n^1-u^1_{0, n}\|_{B_{p, r}^{s-1}}\|\rho^1_{0, n}\|_{B_{p, r}^{s-1}}.
\end{align*}
Denote
\begin{eqnarray*}
X_s=\|\epsilon\|_{B_{p, r}^s}+\|\delta\|_{B_{p, r}^{s-1}},
\end{eqnarray*}
taking the above estimates into (\ref{eq3.6})-(\ref{eq3.7}), we get
\begin{eqnarray*}
X_{s-1}\lesssim \int_0^tX_{s-1}(\|u_n^1,\; u^1_{0, n}\|_{B_{p, r}^s}+\|\rho_n^1,\;\rho^1_{0, n}\|_{B_{p, r}^{s-1}})+t2^{-n(s-\frac 12)},
\end{eqnarray*}
since $\{u_n^1, \rho_n^1\}$ is bounded in $B_{p, r}^s\times B_{p, r}^{s-1}$,
which together with the Gronwall Lemma imply
$$ X_{s-1}\lesssim 2^{-n(s-\frac 12)}.$$
Combining with (\ref{eq3.2}) for $k=1$ and the interpolation inequality, we obtain that
 $$X_s\lesssim X_{s-1}^{\frac{1}{2}}X_{s+1}^{\frac{1}{2}}\lesssim 2^{-\frac{n}{2}(s-\frac 12)}2^{\frac{n}{2}}\lesssim 2^{-\frac{n}{2}(s-\frac 32)}.$$
Thus we have complete the proof of Proposition \ref{pro1}.

In order to obtain the non-uniformly continuous dependence property for the system (\ref{eq4}), we will show that for the constructed  initial data $(u_n^2(0), \rho_n^2(0))$ with small perturbation, it can not approximate to the solution
$(u_n^2, \rho_n^2)$.
\begin{proposition}\label{pro2}
 Under the assumptions of Theorem \ref{the1.1}, we have
      \begin{eqnarray}\label{eq3.9}
      \|u_n^2-u_{0, n}^2-t\mathbf{v}_0^n\|_{B_{p, r}^s}+\|\rho_n^2-\rho_{0, n}^2-t\mathbf{w}_0^n\|_{B_{p, r}^{s-1}}\lesssim t^2+2^{-n\min\{s-\frac 32,\; \frac 12\}},
        \end{eqnarray}
here, $\mathbf{v}_0^n=u_{0, n}^2\partial_xu_{0, n}^2,$ $\mathbf{w}_0^n=k_3u_{0, n}^2\partial_x\rho_{0, n}^2$.
\end{proposition}
\noindent{\bf Proof} \; Firstly, due to (\ref{eq3.1}) and making full use of the product estimates in Lemma \ref{lem2.1}, for $\sigma\geq -1$, we have
\begin{align}
  \|\mathbf{v}_0^n\|_{B_{p, r}^{s+\sigma}}&\lesssim \|u_{0, n}^2\|_{L^\infty}\|\partial_xu_{0, n}^2\|_{B_{p, r}^{s+\sigma}}+ \|u_{0, n}^2\|_{B_{p, r}^{s+\sigma}}\|\partial_xu_{0, n}^2\|_{L^\infty}\nonumber\\
  &\lesssim (2^{-ns}+2^{-n})2^{n(\sigma+1)}+2^{n\sigma}(2^{-ns}2^n+2^{-n})\lesssim 2^{n\sigma}\label{eq3.9}\\
   \|\mathbf{w}_0^n\|_{B_{p, r}^{s+\sigma}}&\lesssim \|u_{0, n}^2\|_{L^\infty}\|\partial_x\rho_{0, n}^2\|_{B_{p, r}^{s+\sigma}}+ \|u_{0, n}^2\|_{B_{p, r}^{s+\sigma}}\|\partial_x\rho_{0, n}^2\|_{L^\infty}\nonumber\\
  &\lesssim (2^{-ns}+2^{-n})2^{n(\sigma+2)}+2^{n\sigma}(2^n2^{-ns}2^n+2^{-n})\lesssim 2^{n(\sigma+1)},\label{eq3.10}\\
\|\mathbf{w}_0^n\|_{B_{p, r}^{s-2}}&\lesssim \|\partial_x\rho_{0, n}^2\|_{B_{p, r}^{s-2}}\|u_{0, n}^2\|_{B_{p, r}^{s-1}}\lesssim 2^{-n}.\label{eq3.11}
\end{align}

Let
\begin{eqnarray*}
        \left\{\begin{array}{ll}
        z_n=u_n^2-u_{0, n}^2-t\mathbf{v}_0^n,\\
         \omega_n=\rho_n^2-\rho_{0, n}^2-t\mathbf{w}_0^n,\end{array}\right.
        \end{eqnarray*}
then we can derive from (\ref{eq4}) that $(z_n, \omega_n)$ satisfies
\begin{eqnarray}\label{eq3.12}
        \left\{\begin{array}{ll}
        \partial_tz_n-u_n^2\partial_xz_n=(z_n+t\mathbf{v}_0^n)\partial_xu_{0, n}^2+tu_n^2\partial_x\mathbf{v}_0^n+f(u_n^2)+g(\rho_n^2)\\
         \partial_t\omega_n-k_3u_n^2\partial_x\omega_n=k_3(z_n+t\mathbf{v}_0^n)\partial_x\rho_{0, n}^2+k_3tu_n^2\partial_x\mathbf{w}_0^n+k_3\rho_n^2\partial_xu_n^2,\\
         z_n(0, x)=0, \; \omega(0, x)=0.
         \end{array}\right.
        \end{eqnarray}
Applying  Lemma \ref{lem2.1}, using (\ref{eq3.3}), (\ref{eq3.9})-(\ref{eq3.10}), we arrive at
\begin{align}
  \|z_n\partial_x\rho_{0, n}^2\|_{B_{p, r}^{s-1}}&\lesssim \|z_n\|_{B_{p, r}^{s-1}}\|\partial_x\rho_{0, n}^2\|_{B_{p, r}^{s-1}}\lesssim 2^n\|z_n\|_{B_{p, r}^{s-1}}, \label{eq3.13}\\
 \|z_n\partial_x\rho_{0, n}^2\|_{B_{p, r}^{s-2}}&\lesssim \|\partial_x\rho_{0, n}^2\|_{B_{p, r}^{s-2}}\|z_n\|_{B_{p, r}^{s-1}}\lesssim \|z_n\|_{B_{p, r}^{s-1}}, \label{eq3.14}\\
  \|\mathbf{v}_0^n\partial_x\rho_{0, n}^2\|_{B_{p, r}^{s-1}}&\lesssim \|\mathbf{v}_0^n\|_{B_{p, r}^{s-1}}\|\partial_x\rho_{0, n}^2\|_{B_{p, r}^{s-1}}\lesssim 2^{-n}2^n\leq C, \label{eq3.15}\\
   \|\mathbf{v}_0^n\partial_x\rho_{0, n}^2\|_{B_{p, r}^{s-2}}&\lesssim \|\partial_x\rho_{0, n}^2\|_{B_{p, r}^{s-2}}\|\mathbf{v}_0^n\|_{B_{p, r}^{s-1}}\lesssim 2^{-n}, \label{eq3.16}\\
  \|u_n^2\partial_x\mathbf{w}_0^n\|_{B_{p, r}^{s-1}}&\lesssim \|u_n^2\|_{B_{p, r}^{s-1}}\|\partial_x\mathbf{w}_0^n\|_{B_{p, r}^{s-1}}\lesssim 2^{-n}2^n\leq C, \label{eq3.17} \\
  \|u_n^2\partial_x\mathbf{w}_0^n\|_{B_{p, r}^{s-2}}&\lesssim \|\partial_x\mathbf{w}_0^n\|_{B_{p, r}^{s-2}}\|u_n^2\|_{B_{p, r}^{s-1}}\lesssim 2^{-n}. \label{eq3.18}
\end{align}
It needs to pay more attention to deal with the term $k_3\rho_n^1\partial_xu_n^1$ and it can be can be decomposed as
\begin{equation*}
 \rho_n^2\partial_xu_n^2=\rho_n^2\partial_xz_n+(\omega_n+\rho_{0, n}^2)\partial_xu_{0, n}^2+t(\rho_n^2\partial_x\mathbf{v}_0^n+\mathbf{w}_0^n\partial_xu_{0 ,n}^2).
\end{equation*}
With (\ref{eq3.3}), (\ref{eq3.9})-(\ref{eq3.10}) at hand, by Lemma \ref{lem2.1}, we find that
\begin{align}
  \|\rho_n^2\partial_xz_n\|_{B_{p, r}^{s-1}}&\lesssim \|\rho_n^2\|_{B_{p, r}^{s-1}}\|\partial_xz_n\|_{B_{p, r}^{s-1}}\lesssim \|z_n\|_{B_{p, r}^s}, \label{eq3.19}\\
  \|\rho_n^2\partial_xz_n\|_{B_{p, r}^{s-2}}&\lesssim \|\partial_xz_n\|_{B_{p, r}^{s-2}}\|\rho_n^2\|_{B_{p, r}^{s-1}}\lesssim \|z_n\|_{B_{p, r}^{s-1}}, \label{eq3.20}\\
  \|\omega_n\partial_xu_{0, n}^2\|_{B_{p, r}^{s-1}}&\lesssim \|\omega_n\|_{B_{p, r}^{s-1}}\|\partial_xu_{0, n}^2\|_{B_{p, r}^{s-1}}\lesssim \|\omega_n\|_{B_{p, r}^{s-1}}, \label{eq3.21}\\
  \|\omega_n\partial_xu_{0, n}^2\|_{B_{p, r}^{s-2}}&\lesssim \|\omega_n\|_{B_{p, r}^{s-2}}\|\partial_xu_{0, n}^2\|_{B_{p, r}^{s-1}}\lesssim \|\omega_n\|_{B_{p, r}^{s-2}}, \label{eq3.22}\\
\|\rho_{0, n}^2\partial_xu_{0, n}^2\|_{B_{p, r}^{s-1}}&\lesssim \|\rho_{0, n}^2\|_{L^\infty}\|\partial_xu_{0, n}^2\|_{B_{p, r}^{s-1}}+ \|\rho_{0, n}^2\|_{B_{p, r}^{s-1}}\|\partial_xu_{0, n}^2\|_{L^\infty}\nonumber\\
  &\lesssim (2^{-ns}2^n+2^{-n})\cdot 1+1\cdot(2^{-ns}2^n+2^{-n})\lesssim 2^{-n\min\{s-1,\;1\}}\label{eq3.23}\\
\|\rho_{0, n}^2\partial_xu_{0, n}^2\|_{B_{p, r}^{s-2}}&\lesssim \|\rho_{0, n}^2\partial_xu_{0, n}^2\|_{B_{p, r}^{s-\frac 32}}\lesssim\|\rho_{0, n}^2\|_{L^\infty}\|\partial_xu_{0, n}^2\|_{B_{p, r}^{s-\frac 32}}+ \|\rho_{0, n}^2\|_{B_{p, r}^{s-\frac 32}}\|\partial_xu_{0, n}^2\|_{L^\infty}\nonumber\\
&\lesssim (2^{-ns}2^n+2^{-n})\cdot 2^{-\frac n2}+2^{-\frac n2}\cdot(2^{-ns}2^n+2^{-n})\lesssim 2^{-n\min\{s-\frac 12,\;\frac 32\}},\label{eq3.24}\\
  \|\rho_n^2\partial_x\mathbf{v}_0^n\|_{B_{p, r}^{s-1}}&\lesssim \|\rho_n^2\|_{B_{p, r}^{s-1}}\|\partial_x\mathbf{v}_0^n\|_{B_{p, r}^{s-1}}\leq C, \label{eq3.25}\\
  \|\rho_n^2\partial_x\mathbf{v}_0^n\|_{B_{p, r}^{s-2}}&\lesssim \|\rho_n^2\|_{B_{p, r}^{s-2}}\|\partial_x\mathbf{v}_0^n\|_{B_{p, r}^{s-1}}\lesssim 2^{-n}, \label{eq3.26}\\
\|\mathbf{w}_0^n\partial_xu_{0 ,n}^2\|_{B_{p, r}^{s-1}}&\lesssim \|\mathbf{w}_0^n\|_{B_{p, r}^{s-1}}\|\partial_xu_{0 ,n}^2\|_{B_{p, r}^{s-1}}\leq C, \label{eq3.27}\\
  \|\mathbf{w}_0^n\partial_xu_{0 ,n}^2\|_{B_{p, r}^{s-2}}&\lesssim \|\mathbf{w}_0^n\|_{B_{p, r}^{s-2}}\|\partial_xu_{0 ,n}^2\|_{B_{p, r}^{s-1}}\lesssim 2^{-n}, \label{eq3.28}
\end{align}

Applying  Lemma \ref{lem2.2} to the second equation of (\ref{eq3.12}), using the fact that $\{u_n^2\}$ is bounded in $L^\infty_T(B_{p, r}^s)$, firstly with (\ref{eq3.14}), (\ref{eq3.16}), (\ref{eq3.18}), (\ref{eq3.20}), (\ref{eq3.22}), (\ref{eq3.24}), (\ref{eq3.26}), (\ref{eq3.28}), we infer that
\begin{equation}\label{eq3.29}
 \|\omega_n\|_{B_{p, r}^{s-2}}\leq C\int_0^t(\|z_n\|_{B_{p, r}^{s-1}}+\|\omega_n\|_{B_{p, r}^{s-2}})d\tau+Ct^22^{-n}+C2^{-n\min\{s-\frac 12,\; \frac 32\}},
\end{equation}
and again combining with (\ref{eq3.13}), (\ref{eq3.15}), (\ref{eq3.17}), (\ref{eq3.19}), (\ref{eq3.21}), (\ref{eq3.23}), (\ref{eq3.25}), (\ref{eq3.27}), we obtain that
\begin{equation}\label{eq3.30}
 \|\omega_n\|_{B_{p, r}^{s-1}}\leq C\int_0^t(\|z_n\|_{B_{p, r}^s}+\|\omega_n\|_{B_{p, r}^{s-1}})d\tau+C\int_0^t2^n\|z_n\|_{B_{p, r}^{s-1}}d\tau+Ct^2+C2^{-n\min\{s-1,\; 1\}}.
\end{equation}

In the following, we shall estimate $z_n$ in $B_{p, r}^{s-1}$ and $B_{p, r}^s$, respectively.
With the aid of Lemma \ref{lem2.1} and (\ref{eq3.3}), (\ref{eq3.9}), one has
\begin{align}
 \|z_n\partial_xu_{0, n}^2\|_{B_{p, r}^s}&\lesssim \|z_n\|_{L^\infty}\|\partial_xu_{0, n}^2\|_{B_{p, r}^s}
+\|z_n\|_{B_{p, r}^s}\|\partial_xu_{0, n}^2\|_{L^\infty}\nonumber\\
&\lesssim \|z_n\|_{B_{p, r}^{s-1}}\|\partial_xu_{0, n}^2\|_{B_{p, r}^s}
+\|z_n\|_{B_{p, r}^s}\|\partial_xu_{0, n}^2\|_{B_{p, r}^{s-1}}\nonumber\\
&\lesssim 2^n\|z_n\|_{B_{p, r}^{s-1}}+\|z_n\|_{B_{p, r}^s},\label{eq3.31}\\
 \|z_n\partial_xu_{0, n}^2\|_{B_{p, r}^{s-1}}&\lesssim \|z_n\|_{B_{p, r}^{s-1}}\|\partial_xu_{0, n}^2\|_{B_{p, r}^{s-1}}\lesssim \|z_n\|_{B_{p, r}^{s-1}},\label{eq3.32}\\
\|\mathbf{v}_0^n\partial_xu_{0, n}^2\|_{B_{p, r}^s}
&\lesssim \|\mathbf{v}_0^n\|_{B_{p, r}^{s-1}}\|\partial_xu_{0, n}^2\|_{B_{p, r}^s}
+\|\mathbf{v}_0^n\|_{B_{p, r}^s}\|\partial_xu_{0, n}^2\|_{B_{p, r}^{s-1}}\leq C,\label{eq3.33}\\
\|\mathbf{v}_0^n\partial_xu_{0, n}^2\|_{B_{p, r}^{s-1}}
&\lesssim \|\mathbf{v}_0^n\|_{B_{p, r}^{s-1}}\|\partial_xu_{0, n}^2\|_{B_{p, r}^{s-1}}\lesssim 2^{-n},\label{eq3.34}\\
\|u_{0, n}^2\partial_x\mathbf{v}_0^n\|_{B_{p, r}^s}
&\lesssim \|u_{0, n}^2\|_{B_{p, r}^{s-1}}\|\partial_x\mathbf{v}_0^n\|_{B_{p, r}^s}
+\|u_{0, n}^2\|_{B_{p, r}^s}\|\partial_x\mathbf{v}_0^n\|_{B_{p, r}^{s-1}}\leq C,\label{eq3.35}\\
\|u_{0, n}^2\partial_x\mathbf{v}_0^n\|_{B_{p, r}^{s-1}}&\lesssim\|u_{0, n}^2\|_{B_{p, r}^{s-1}}\|\partial_x\mathbf{v}_0^n\|_{B_{p, r}^{s-1}}\lesssim 2^{-n}. \label{eq3.36}
\end{align}

For the term $f(u_n^2)=f_1(u_n^2)+f_2(u_n^2),$  we have from Lemma \ref{lem2.1} and (\ref{eq3.3}) that
  \begin{equation}\label{eq3.37}
    \|f_1(u_n^2)\|_{B_{p, r}^s}\lesssim \|(u_n^2)^2\|_{B_{p, r}^{s-1}}\lesssim  \|u_n^2\|^2_{B_{p, r}^{s-1}}\lesssim 2^{-2n},
  \end{equation}
while it needs to be more careful to deal with $f_2(u_n^2)$. By making full use of the structure of $u_n^2$, we find that
\begin{align*}
f_2(u_n^2)=&\underbrace{\frac {3-k_1}{2}\partial_x(1-\partial_x^2)^{-1}(\partial_x(u_n^2+u_{0, n}^2)\partial_xz_n)}_{f_{2, 1}}\\
&\;+\underbrace{\frac {3-k_1}{2}\partial_x(1-\partial_x^2)^{-1}(t\partial_x(u_n^2+u_{0, n}^2)\partial_x\mathbf{v}_0^n)}_{f_{2, 2}}\\
&\;+\underbrace{\frac {3-k_1}{2}\partial_x(1-\partial_x^2)^{-1}((\partial_xu^2_{0, n})^2)}_{f_{2, 3}},
\end{align*}
and
\begin{align}
  \|f_{2, 1}\|_{B_{p, r}^s}&\lesssim \|\partial_x(u_n^2+u_{0, n}^2)\partial_xz_n\|_{B_{p, r}^{s-1}}\lesssim \|z_n\|_{B_{p, r}^s}, \label{eq3.38}\\
  \|f_{2, 1}\|_{B_{p, r}^{s-1}}&\lesssim \|\partial_x(u_n^2+u_{0, n}^2)\partial_xz_n\|_{B_{p, r}^{s-2}}\lesssim \|z_n\|_{B_{p, r}^{s-1}}, \label{eq3.39}\\
 \|f_{2, 2}\|_{B_{p, r}^s}&\lesssim \|\partial_x(u_n^2+u_{0, n}^2)\partial_x\mathbf{v}_0^n\|_{B_{p, r}^{s-1}}\leq C,\label{eq3.40}\\
  \|f_{2, 2}\|_{B_{p, r}^{s-1}}&\lesssim \|\partial_x(u_n^2+u_{0, n}^2)\partial_x\mathbf{v}_0^n\|_{B_{p, r}^{s-2}}\lesssim 2^{-n},\label{eq3.41}\\
   \|f_{2, 3}\|_{B_{p, r}^s}&\lesssim  \|(\partial_xu^2_{0, n})^2\|_{B_{p, r}^{s-1}}\lesssim\|\partial_xu^2_{0, n}\|_{L^\infty}\|\partial_xu^2_{0, n}\|_{B_{p, r}^{s-1}}\nonumber\\
   &\lesssim(2^{-ns}2^n+2^{-n})\cdot 1\lesssim 2^{-n\min\{s-1, \; 1\}}, \label{eq3.42}\\
\|f_{2, 3}\|_{B_{p, r}^{s-1}}&\lesssim  \|(\partial_xu^2_{0, n})^2\|_{B_{p, r}^{s-\frac 32}}\lesssim\|\partial_xu^2_{0, n}\|_{L^\infty}\|\partial_xu^2_{0, n}\|_{B_{p, r}^{s-\frac 32}}\nonumber\\
   &\lesssim(2^{-ns}2^n+2^{-n})\cdot 2^{-\frac n2}\lesssim 2^{-n\min\{s-\frac 12, \; \frac 32\}}. \label{eq3.43}
\end{align}
$g(\rho_n^2)$ can be performed in a similar way. Firstly, we obtain that
\begin{align*}
g(\rho_n^2)=&\underbrace{\frac {k_2}{2}\partial_x(1-\partial_x^2)^{-1}(\omega_n(\rho_n^2+\rho_{0, n}^2))}_{g_1}+\underbrace{\frac {k_2}{2}\partial_x(1-\partial_x^2)^{-1}(t\mathbf{w}_0^n(\rho_n^2+\rho_{0, n}^2))}_{g_2}\\
&\;+\underbrace{\frac {k_2}{2}\partial_x(1-\partial_x^2)^{-1}((\rho^2_{0, n})^2)}_{g_3},
\end{align*}
and
\begin{align}
 \|g_1\|_{B_{p, r}^s}&\lesssim \|\omega_n(\rho_n^2+\rho_{0, n}^2)\|_{B_{p, r}^{s-1}}\lesssim \|\omega_n\|_{B_{p, r}^{s-1}}, \label{eq3.44}\\
  \|g_1\|_{B_{p, r}^{s-1}}&\lesssim \|\omega_n(\rho_n^2+\rho_{0, n}^2)\|_{B_{p, r}^{s-2}}\lesssim \|\omega_n\|_{B_{p, r}^{s-2}}, \label{eq3.45}\\
 \|g_2\|_{B_{p, r}^s}&\lesssim \|\mathbf{w}_0^n(\rho_n^2+\rho_{0, n}^2)\|_{B_{p, r}^{s-1}}\leq C,\label{eq3.46}\\
  \|g_2\|_{B_{p, r}^{s-1}}&\lesssim \|\mathbf{w}_0^n(\rho_n^2+\rho_{0, n}^2)\|_{B_{p, r}^{s-2}}\lesssim 2^{-n},\label{eq3.47}\\
   \|g_3\|_{B_{p, r}^s}&\lesssim  \|(\rho^2_{0, n})^2\|_{B_{p, r}^{s-1}}\lesssim\|\rho^2_{0, n}\|_{L^\infty}\|\rho^2_{0, n}\|_{B_{p, r}^{s-1}}\nonumber\\
   &\lesssim(2^{-ns}2^n+2^{-n})\cdot 1\lesssim 2^{-n\min\{s-1, \; 1\}}, \label{eq3.48}
   \end{align}
   \begin{align}
\|g_3\|_{B_{p, r}^{s-1}}&\lesssim  \|(\rho^2_{0, n})^2\|_{B_{p, r}^{s-\frac 32}}\lesssim\|\rho u^2_{0, n}\|_{L^\infty}\|\rho^2_{0, n}\|_{B_{p, r}^{s-\frac 32}}\nonumber\\
   &\lesssim(2^{-ns}2^n+2^{-n})\cdot 2^{-\frac n2}\lesssim 2^{-n\min\{s-\frac 12, \; \frac 32\}}. \label{eq3.49}
\end{align}
Applying Lemma \ref{lem2.2} firstly together with (\ref{eq3.32}), (\ref{eq3.34}), (\ref{eq3.36}),  (\ref{eq3.37}), (\ref{eq3.39}), (\ref{eq3.41}), (\ref{eq3.43}),  (\ref{eq3.45}), (\ref{eq3.47}), (\ref{eq3.49}) to the first equation of (\ref{eq3.12}), using the fact that $\{u^2_n\}$ is bounded in $L^\infty_T(B_{p, r}^s)$, we infer that
\begin{equation}\label{eq3.50}
 \|z_n\|_{B_{p, r}^{s-1}}\leq C\int_0^t(\|z_n\|_{B_{p, r}^{s-1}}+\|\omega_n\|_{B_{p, r}^{s-2}})d\tau+Ct^22^{-n}+C2^{-n\min\{s-\frac 12,\; \frac 32\}},
\end{equation}
and again combining with (\ref{eq3.31}), (\ref{eq3.33}), (\ref{eq3.35}),  (\ref{eq3.37}), (\ref{eq3.38}), (\ref{eq3.40}), (\ref{eq3.42}),  (\ref{eq3.44}), (\ref{eq3.46}), (\ref{eq3.48}), we obtain
\begin{equation}\label{eq3.51}
 \|z_n\|_{B_{p, r}^s}\leq C\int_0^t(\|z_n\|_{B_{p, r}^s}+\|\omega_n\|_{B_{p, r}^{s-1}})d\tau+ C\int_0^t2^{n}\|z_n\|_{B_{p, r}^{s-1}}d\tau+Ct^2+C2^{-n\min\{s-1, \; 1\}}.
\end{equation}
Using Gronwall Lemma to (\ref{eq3.23}) and (\ref{eq3.43}) imply
\begin{align*}
 \|z_n\|_{B_{p, r}^{s-1}}+\|\omega_n\|_{B_{p, r}^{s-2}}\leq Ct^22^{-n}+C2^{-n\min\{s-\frac 12,\; \frac 32\}},
\end{align*}
which together with (\ref{eq3.30}) and (\ref{eq3.51}) yield that
 \begin{equation*}
  \|z_n\|_{B_{p, r}^s}+\|\omega_n\|_{B_{p, r}^{s-1}}\leq Ct^2+C2^{-n\min\{s-\frac 32,\; \frac 12\}}.
\end{equation*}
 Thus, we have finished the proof of Proposition \ref{pro2}.

{\bf Proof of Theorem \ref{the1.1}} It is obvious that
 \begin{align*}
   \|u_{0, n}^2-u_{0, n}^1\|_{B_{p, r}^s}&=\|g_n\|_{B_{p, r}^s}\leq C2^{-n},\\
   \|\rho_{0, n}^2-\rho_{0, n}^1\|_{B_{p, r}^{s-1}}&=\|g_n\|_{B_{p, r}^{s-1}}\leq C2^{-n},
 \end{align*}
 which means that
 \begin{equation*}
   \lim_{n\rightarrow \infty}(\|u_{0, n}^2-u_{0, n}^1\|_{B_{p, r}^s}+ \|\rho_{0, n}^2-\rho_{0, n}^1\|_{B_{p, r}^{s-1}})=0.
 \end{equation*}
However, according to Proposition \ref{pro1} and Proposition \ref{pro2}, we get
\begin{align}\label{eq3.52}
  \|\rho_n^2-\rho_n^1\|_{B_{p, r}^{s-1}}&=\|\omega_n+t\mathbf{w}_0^n+g_n+\rho_{0, n}^1-\rho_n^1\|_{B_{p, r}^{s-1}}\nonumber\\
  &\gtrsim t\|\mathbf{w}_0^n\|_{B_{p, r}^{s-1}}-2^{-n}-t^2-2^{-n\min\{s-\frac 32,\; \frac 12\}}-2^{-\frac n2(s-\frac 32)}.
\end{align}
Notice that
\begin{align*}
 \mathbf{w}_0^n&=k_3u_{0, n}^2\partial_x\rho_{0, n}^2=k_3(f_n+g_n)\partial_x(2^nf_n+g_n)\\
 &=k_3f_n\partial_x(2^nf_n)+k_3g_n\partial_x(2^nf_n)+k_3f_n\partial_xg_n+k_3g_n\partial_xg_n.
\end{align*}
With the aid of Lemma \ref{lem2.1} and the Banach algebra property of $B_{p, r}^{s-1}$, we find that
\begin{align*}
 \|f_n\partial_x(2^nf_n)\|_{B_{p, r}^{s-1}}&\lesssim \|f_n\|_{L^\infty}\|\partial_x(2^nf_n)\|_{B_{p, r}^{s-1}}+\|f_n\|_{B_{p, r}^{s-1}}\|\partial_x(2^nf_n)\|_{L^\infty}\\
 &\lesssim 2^{-ns}2^n+2^{-n}2^n2^{-ns}2^n\lesssim 2^{-n(s-1)},\\
\|f_n\partial_xg_n\|_{B_{p, r}^{s-1}}&\lesssim \|f_n\|_{B_{p, r}^{s-1}}\|\partial_xg_n\|_{B_{p, r}^{s-1}}\lesssim 2^{-2n},\\
\|g_n\partial_xg_n\|_{B_{p, r}^{s-1}}&\lesssim \|g_n\|_{B_{p, r}^{s-1}}\|\partial_xg_n\|_{B_{p, r}^{s-1}}\lesssim 2^{-2n}.
\end{align*}
However, using the fact that $\Delta_j\big(g_n\partial_x(2^nf_n)\big)=0, j\neq n$ and $\Delta_n\big(g_n\partial_x(2^nf_n)\big)=g_n\partial_x(2^nf_n)$ for $n\geq 5$, direct calculation
 shows that for $n\gg 1$,
\begin{align*}
  &\|g_n\partial_x(2^nf_n)\|_{B_{p, r}^{s-1}}=2^{n(s-1)}\|g_n\partial_x(2^nf_n)\|_{L^p}\\
  =&\|2^{-n}\phi(x)\partial_x\phi(x)\sin(\frac{17}{12}2^nx)+\frac{17}{12}\phi^2(x)\cos(\frac{17}{12}2^nx)\|_{L^p}\\
  \gtrsim &\|\frac{17}{12}\phi^2(x)\cos(\frac{17}{12}2^nx)\|_{L^p}-2^{-n}\rightarrow \frac{17}{12} \big(\frac{\int_0^{2\pi}|\cos x|^pdx}{2\pi}\big)^{\frac 1p}\|\phi^2(x)\|_{L^p},
\end{align*}
by the Riemann Theorem.

 Taking the above estimates into (\ref{eq3.52}) yields
  \begin{align*}
   \liminf_{n\rightarrow \infty}\|\rho_n^2-\rho_n^1\|_{B^{s-1}_{p,r}}\gtrsim t\quad\text{for} \ t \ \text{small enough}.
  \end{align*}
 Similarly, we have
   \begin{align*}
   \liminf_{n\rightarrow \infty}\|u_n^2-u_n^1\|_{B^s_{p,r}}\gtrsim t\quad\text{for} \ t \ \text{small enough}.
  \end{align*}

  This completes the proof of Theorem \ref{the1.1}.

\section*{Acknowledgments}
 This work is supported by the National Natural Science Foundation of China (Grant No.12001163).

\end{document}